\documentclass[11pt,a4paper]{amsart}

\usepackage{graphicx}

\newtheorem{theorem}{Theorem}[section]

\newtheorem{prop}[theorem]{Proposition}

\theoremstyle{definition}

\theoremstyle{remark}
\newtheorem{remark}[theorem]{Remark}

\renewcommand{\theequation}{\arabic{section}.\arabic{equation}}

\def\Q{{\mathbb{Q}}}  \def\C{{\mathbb{C}}}
\def\Z{{\mathbb{Z}}}  \def\P{{\mathbb{P}}}
\def\A{{\mathbb{A}}} \def\F{{\mathbb{F}}}

\newcommand{\rank}{\operatorname{rank}}

\def\a{\alpha} \def\b{\beta} \def\g{\gamma}  
\def\io{\iota}   
   \def\s{\sigma} \def\t{\tau}
  \def\w{\omega}  

\def\phihat{\kern.3em\hat{\kern-.3em\phi}}
\def\sphihat{\kern.2em\hat{\kern-.2em\phi}}
\def\e{\varepsilon}
\def\ph{\varphi}

\def\<{\langle} \def\>{\rangle}

\def\ds{\displaystyle}

\def\tS{S}
\def\tr{\mbox{\scriptsize tr}}
\newcommand{\leg}[2]{\left(\frac{#1}{#2}\right)}

\makeatletter
\def\currenteq{\textup{\tagform@{\p@equation\theequation} }}
\makeatother

\title
{A singular $K3$ surface related to \\ sums of consecutive cubes}

\author{Masato Kuwata}
\address{D\'epartement de Math\'ematiques\\
Universit\'e de Caen\\
B.P. 5186\\
F-14032 Caen Cedex\\
France}
\email{kuwata@math.unicaen.fr}

\author{Jaap Top}
\address{Vakgroep Wiskunde R{\sl u}G \\
P.O. Box 800 \\
9700 AV Groningen \\
the Netherlands}

\email{top@math.rug.nl}

\date{December 1999, revised version}

\begin{document}

\maketitle

\section{Introduction}
The well known formula
$$1^3+2^3+\ldots+k^3=(k(k+1)/2)^2$$
shows in particular that the sum of the first $k$
consecutive cubes is in fact a square. A lot of related
diophantine problems have been studied. For instance,
Stroeker \cite{Str} considered the question, for which
integers $m>1$ one can find non-trivial solutions
of $m^3+(m+1)^3+\ldots+(m+k-1)^3=\ell^2$. In particular
he found all solutions with $1<m<51$ and with $m=98$.
Several authors considered the question of finding sums of
consecutive squares which are themselves squares;
see, e.g., \cite{P-R} and references given there.
In our paper \cite{K-T} we considered this problem
from a geometric point of view by exploiting a connection
between solutions and rational points on the rational
elliptic surface given by $6y^2=(x^3-x)-(t^3-t)$.

Here we present a similar treatment of the
equation
\begin{equation}\label{eq_mkl}
m^{3} + (m+1)^{3} + (m+2)^{3} + \dots + (m+k-1)^{3} =l^3.
\end{equation}

It turns out that \currenteq determines a so-called
singular $K3$ surface.  Although our emphasis will be on
geometric and arithmetic algebraic properties of this
surface, we will note some consequences for the
diophantine equation. In particular, we show that infinitely many
non-trivial solutions exist; a result that was already known
to C.~Pagliani in 1829/30. Pagliani's solutions can be
regarded as points on a rational curve in the $K3$ surface. 

By putting
\begin{equation}\label{chg_var}
x=k,\qquad
y=2m+k-1,\qquad
z=2l,
\end{equation}
the equation \eqref{eq_mkl}  becomes
\begin{equation}\label{eq_xyz}
xy(x^2+y^2-1)=z^3.
\end{equation}
This equation \eqref{eq_xyz} determines an algebraic surface in
the affine space $\A^3$.  We consider the projective model
of it, given by
\[
S_0:\bigl\{(X:Y:Z:W)\mid XY(X^2+Y^2-W^2)=Z^3W\bigr\} \subset
\P^3.
\]
This surface has some singular points.  Denote by
$S$ the minimal non-singular model of $S_0$.  We will show that
the surface $S$ is a $K3$ surface whose Picard number is 20,
which is maximal for a $K3$ surface defined over a field of
characteristic 0. Such $K3$ surfaces with maximal Picard
number have been classified by Shioda and Inose \cite{shioda-inose}.

In sections 2 and 3 we study geometric properties of the surface $S$.
We construct some elliptic fibrations on it and determine the
N\'eron-Severi group of $S$. Moreover, we find the Shioda-Inose class 
of this surface. It turns out that he elliptic fibrations we construct,
do not have sections which would produce infinitely many
non-trivial solutions to the original diophantine equation.
In section 4 we describe a base change of one of the elliptic pencils.
Using this, we do find infinitely many
non-trivial solutions to \eqref{eq_mkl}. It turns out that precisely
the same solutions already appeared in 1829/30 in work of C.~Pagliani.

Section 5 gives a description of the surface $S$ and the
base change we present, in terms of the product of two curves.
This is used in section 6 to determine the Hasse-Weil zeta function 
of $S$ over $\Q$.

{\small
\begin{proof}[Acknowledgements]
The first-named author would like to thank Joe Silverman for useful
suggestions in regard to Proposition~\ref{no-poly-solution}.  He also
thanks the Centre Interunivesitaire en Calcul Math\'ematique
Alg\'ebrique at Concordia University for allowing him to use their
computing facilities. We are grateful to the referee for pointing out
the reference to Pagliani, and for some interesting
remarks concerning the product of two curves we considered.
Finally we thank Ken Ono for mentioning
the formula for our cusp form in terms of eta-functions.
\renewcommand{\qed}{}\end{proof}}

\section{Elementary properties of $S_{0}$}

\subsection{Symmetry}\label{symm}
The original diophantine problem asks for solutions to
\eqref{eq_mkl} in {\em positive} integers.  This restriction,
however, is not an essential one.  For example,
$(m,k,l)=(-2,8,6)$ is a solution to~\eqref{eq_mkl}.  Then,
observing
\[
\underbrace{(-2)^3 + (-1)^3 + 0^3 + 1^3 + 2^3}_{\ds 0}
+3^3+4^3+5^3=6^3,
\]
we find another solution $(m,k,l)=(3,3,6)$
corresponding to $3^3+4^3+5^3=6^3$.  This reflects an instance
of the action of a certain group of
symmetries on the set of solutions; one can always find a
nonnegative solution in the orbit of a solution.
These symmetries are conveniently described in terms of the coordinates
introduced in \eqref{chg_var}. There are obvious involutions acting on
\eqref{eq_xyz}, namely
\begin{align*}
\t_1 &: (x,y,z)\longmapsto(-x,y,-z), \\
\t_2 &: (x,y,z)\longmapsto(x,-y,-z),\\
\t_3 &: (x,y,z)\longmapsto(y,x,z).
\end{align*}
These involutions generate a group of order~$8$, which is isomorphic
to the dihedral group $D_{4}$.  The symmetry may be seen clearly in
Figure~1.

\begin{figure}[h]
\centerline{\includegraphics[scale=.707]{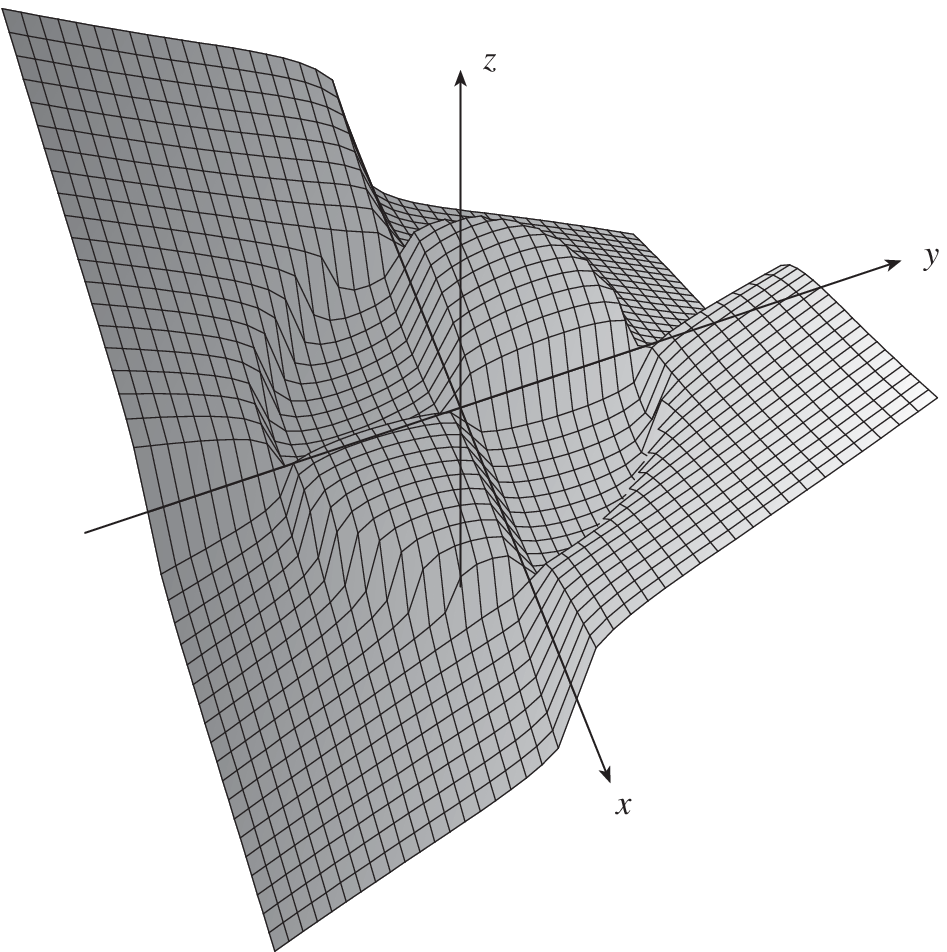}}
\caption{}
\end{figure}

It is easy to see that $(m,k,l)$ is an integral solution to \eqref{eq_mkl}
if and only if $(x,y,z)$ is an integral solution to \eqref{eq_xyz}
such that $x$ and $y$ have different parity.

The surface $S_{0}$ admits another automorphism given by
\[
\t_{4}:(x,y,z) \mapsto (x,y,\w z),
\]
where $\w$ is a primitive third root of unity.
We see that $\t_{4}$ commutes with $\t_{1}$, $\t_{2}$ and $\t_{3}$,
and thus $\t_{1},\dots,\t_{4}$ forms a group isomorphic to
$D_{4}\times C_{3}$, where $C_{3}$ is the cyclic group of order~$3$.
It is not difficult to see that these are the only automorphisms of
$S_{0}$ induced from those of $\P^{3}$.

\subsection{Singularities}
The surface $S_{0}$ has the following five singular points
\[
(X:Y:Z:W)=(0:0:0:1), \ (\pm1:0:0:1), \ (0:\pm1:0:1),
\]
all of whom are rational double points. These are the only
singularities of $S_0$. Hence by the following result 
one concludes that the minimal non-singular model $S$ of $S_0$ is 
a $K3$ surface.

\begin{prop}{\cite[Prop.~2.1]{P-T-vdV}}
Any surface which is the minimal resolution of singularities of a
surface given by a degree four equation in projective 3-space with at
most rational double points as singularities, is a $K3$ surface. \qed
\end{prop}

Note that all of the singularities of $S_{0}$ are of type $A_{2}$,
and thus the fibre of $S\rightarrow S_0$ above each singular point
consists of two rational curves. This can be seen directly, but it 
also follows from the results in the next section since the resolution 
of singularities
is built-in in Tate's algorithm for determining the singular fibers of
an elliptic pencil.

\subsection{Lines contained in $S_{0}$}
The surface $S_{0}$ contains the following lines:
\begin{align*}
\ell_{1} : &\left\{\begin{aligned} X&=0 \\ Z&=0 \end{aligned}\right. &
\ell_{2} : &\left\{\begin{aligned} Y&=0 \\ Z&=0 \end{aligned}\right. &
\ell_{3} : &\left\{\begin{aligned} X&=0 \\ W&=0 \end{aligned}\right. \\
\ell_{4} : &\left\{\begin{aligned} Y&=0 \\ W&=0 \end{aligned}\right. &
\ell_{5} : &\left\{\begin{aligned} X&=iY  \\ W&=0 \end{aligned}\right. &
\ell_{6} : &\left\{\begin{aligned} X&=-iY \\ W&=0 \end{aligned}\right.
\\
\ell_{7} : &\left\{\begin{aligned} X&=W \\ Y&=Z \end{aligned}\right. &
\ell_{8} : &\left\{\begin{aligned} X&=W \\ Y&=\w Z \end{aligned}\right.
&
\ell_{9} : &\left\{\begin{aligned} X&=W \\ Y&=\w^{2}Z
\end{aligned}\right.\\
\ell_{10} : &\left\{\begin{aligned} X&=-W \\ Y&=-Z \end{aligned}\right.
&
\ell_{11} : &\left\{\begin{aligned} X&=-W \\ Y&=-\w Z
\end{aligned}\right. &
\ell_{12} : &\left\{\begin{aligned} X&=-W \\ Y&=-\w^{2}Z
\end{aligned}\right.\\
\ell_{13} : &\left\{\begin{aligned} Y&=W \\ X&=Z \end{aligned}\right. &
\ell_{14} : &\left\{\begin{aligned} Y&=W \\ X&=\w Z \end{aligned}\right.
&
\ell_{15} : &\left\{\begin{aligned} Y&=W \\ X&=\w^{2}Z
\end{aligned}\right.\\
\ell_{16} : &\left\{\begin{aligned} Y&=-W \\ X&=-Z \end{aligned}\right.
&
\ell_{17} : &\left\{\begin{aligned} Y&=-W \\ X&=-\w Z
\end{aligned}\right. &
\ell_{18} : &\left\{\begin{aligned} Y&=-W \\ X&=-\w^{2}Z.
\end{aligned}\right.\\
\end{align*}

The group of symmetries generated by $\{\t_{1},\dots,\t_{4}\}$ acts on
 this set of lines.  It divides the set into  four orbits: 
 $\{\ell_{1},\ell_{2}\}$,
 $\{\ell_{3},\ell_{4}\}$, $\{\ell_{5},\ell_{6}\}$ and
 $\{\ell_{7},\dots,\ell_{18}\}$.
 
With the aid of the description of the N\'{e}ron-Severi group of $S$
in terms of an elliptic fibration, one can in fact show
that these are the only lines contained in $S_0$. However,
we will not need this in the sequel.

\section{Elliptic fibrations on $S$ and the N\'eron-Severi group}

By regarding $y$ as a constant in \eqref{eq_xyz}, one obtains
a plane cubic curve in the $xz$-plane.  It is easy to see that
these cubic curves are nonsingular except for finitely many
values of $y$. 
This implies that the map $S_{0}\to \A^{1}$ given by $(x,y,z)\mapsto y$
determines an elliptic fibration $\e_1: S\to \P^1$.  Moreover, since
$(x,z)=(0,0)$ is a rational point on each fiber, $\e_1$
admits a section $\s_{0}:\P^{1}\to S$, which we designate as the
0-section.

Let $t$ be the generic point of the base curve $\P^{1}$, and let $E_{t}$
be the fiber of $\e_{1}$ at $t$.  This is nothing but the curve defined
over $\Q(t)$ obtained by replacing $y$ by $t$ in \eqref{eq_xyz}.
Using the change of variables
\[
\left\{
\begin{aligned}
x_{1}&=\frac{(t^3-t)z}{x}\\
y_{1}&=\frac{(t^3-t)^2}{x},
\end{aligned}\right.
\]
the curve $E_{t}$ can be written in the Weierstrass form
\begin{equation}\label{weier}
E_t : \ y_{1}^2=x_{1}^3-t^4(t^2-1)^3.
\end{equation}

The bad fibers of the fibration $e_1$ are easily determined
using Tate's algorithm. They are given in the table below. Here, we
denote by $m_t$ (resp. $m_t^{(1)}$)
the number of irreducible (resp. simple) components in the
fiber of $E_t$ at $t$, and we denote by $\chi$ the
(topological) Euler number of the bad fiber.
\[
\begin{tabular}{|c|ccccc|}
\hline
$t$ & Kodaira type & $\chi$ &
$m_t$&
$m_t^{(1)}$&
$\begin{array}{l}
\hbox{Group} \\ \hbox{structure}
\end{array}$ \\[4pt]
\hline
0 & IV$^*$ & 8 & 7 & 3 &${\mathbb{G}}_a \times \Z/3\Z$
\vrule width 0pt height 15pt\\[3pt]
$\pm1$ & I$_0^*$ & 6 & 5 & 4& ${\mathbb{G}}_a \times (\Z/2\Z)^2$
\\[3pt]
$\infty$ & IV & 4 & 3 & 3& ${\mathbb{G}}_a \times \Z/3\Z $
\\[3pt]
\hline \end{tabular}
\]

Here group structure means the structure of a special fiber
of the N\'{e}ron model of $E_t/\C(t)$. 
In order to determine the N\'eron-Severi group of
$S_{\C}=S\times_{\Q}\C$, we need to determine the Mordell-Weil group
$E_t(\C(t))$.  We observe that the image of the lines
$\ell_{1},\dots,\ell_{6}$ in $S$ are
components of  bad fibers.  It turns out that the lines
$\ell_{7},\dots,\ell_{18}$ determine sections of $\e_{1}$.
For example, $\ell_{7}$ is transformed to the section
\[
\s_1 =\bigl(t^2(t^2-1), t^2(t^2-1)^2\bigr).
\]
Note that this corresponds to the trivial parametric solution
$(m,k,l)=(m,0,m)$ to the original equation~\eqref{eq_mkl}.  Since the
curve $E_t/\C(t)$ has complex multiplication by $\Q(\w)$, we see that
\[
[\w]\s_1 =\bigl(\w t^2(t^2-1),t^2(t^2-1)^2\bigr).
\]
is also a section. It corresponds to
the line $\ell_{8}$ in $S$.

\begin{prop}\label{MWL}
\begin{enumerate}
\item
The Mordell-Weil group $E_t(\C(t))$ is isomorphic to $\Z\oplus\Z$, and
it is generated by $\s_1$ and $[\w]\s_1$.  
\item The Mordell-Weil lattice
of $E_{t}$ over $\C(t)$  is isomorphic to $A_{2}^{*}$, the dual
lattice of the root lattice $A_2$.
\item
The surface $S$ is a singular $K3$ surface.  The determinant of
its N\'eron-Severi lattice is $-48$.
\end{enumerate}
\end{prop}

Recall that the Mordell-Weil lattice in
the sense of Shioda \cite{shioda:MWL} of an elliptic surface
with section consists of the Mordell-Weil group modulo torsion
of its generic fiber, together with twice the canonical height
pairing on it. Also, a $K3$ surface in characteristic zero
is called {\em singular} (or
{\em exceptional\/}) if its N\'eron-Severi group has
the maximal rank, which is $20$.

\begin{proof}  We first determine the torsion
subgroup $E_t(\C(t))_{\hbox{\tiny tors}}$.  It is known
that the specialization map 
$E_t^0(\C(t))_{\hbox{\tiny tors}}\to E_{t_0}^{\hbox{\tiny ns}}(\C)$ 
for any $t_0\in\C$ is always
injective even when the fiber is a bad fiber (see
\cite[Lemma~1.1 (b)]{miranda-persson}); here $E_t^0(\C(t))$
denotes the subgroup of $E_t(\C(t))$ consisting of all points
which specialize to smooth points in $E_{t_0}(\C)$.
Moreover, the notation $E_{t_0}^{\hbox{\tiny ns}}(\C)$
is used for the group of smooth points in $E_{t_0}(\C)$,
which is also the connected component of zero of the fiber
over $t_0$ in the N\'{e}ron model of $E_t$. 
Considering the table of bad
fibers above, one observes that $E_t(\C(t))_{\hbox{\tiny tors}}$ is a
subgroup of both $\C\times \Z/3\Z$ and
$\C\times (\Z/2\Z)^2$.  This implies that
$E_t(\C(t))_{\hbox{\tiny tors}}$ is trivial.

Since $\s_1$ is not the zero section, we conclude that
it has infinite order. Also, since the Mordell-Weil group
is in fact a module over $\mbox{End}(E_t)\cong\Z[\omega]$,
it follows that $\s_1$ and $[\w]\s_1$ are independent.

The rank of the N\'{e}ron-Severi group of $S$ and the rank
of the Mordell-Weil group are related by the Shioda-Tate formula,
which says that
\[
\rank NS(S) = 2 + \sum_t (m_t-1) + \rank E_t(\C(t)).
\]
{}From the calculation of the bad fibers and the fact that
$\rank NS(S)\leq 20$ one concludes that the
rank of $E_t(\C(t))$ is at most $2$, hence equal to $2$.
It follows that the rank of $NS(S)$ is~20 and therefore
$S$ is a singular $K3$ surface. 

Next we compute the canonical height pairing for $\{\s_1,
[\w]\s_1\}$.  This is easily done using \cite{kuwata:can-height}, and
one finds the matrix
\[\def\arraystretch{1.1}
\left(\begin{array}{rr}
\frac{1}{3} & -\frac{1}{6}\\
-\frac{1}{6} & \frac{1}{3}
\end{array}\right).
\]
This shows again that $\s_1$ and $[\w]\s_1$ are linearly
independent. We use the computation to conclude that
$\{\s_1,[\w]\s_1\}$ generates the
Mordell-Weil group. If the Mordell-Weil group $E_t(\C(t))$
would be  strictly larger than the group generated by $\s_1$ and
$[\w]\s_1$, then $E_t(\C(t))$ must contain an element whose
canonical height is less than $1/3$.  However, the a priori
lower bound for the canonical height calculated by using the
method of \cite{kuwata:can-height} is $1/3$ in the present case. 
Hence $\s_1$ and $[\w]\s_1$ generate $E_t(\C(t))$.

Let us denote by $MW$ the Mordell-Weil lattice
of $E_t$ in the sense of Shioda\cite{shioda:MWL}.  The pairing on
$MW$ is given by twice the canonical height pairing; i.e.,
it is given by the matrix
\[\def\arraystretch{1.2}
\left(\begin{array}{rr}
\frac{2}{3} & -\frac{1}{3}\\
-\frac{1}{3} & \frac{2}{3}
\end{array}\right).
\]
Hence, $MW$ is isomorphic to $A_{2}^{*}$, and we have
\[
|\det NS(S)| =\frac{\prod m_t^{(1)} |\det MW|}
{|E_t(\C(t))_{\hbox{\tiny tors}}|^2}
=\frac{3\cdot 4 \cdot 4 \cdot 3 \cdot \frac{1}{3}}{1}
=48.
\]
Since by the Hodge index theorem $\det NS(S)$ is negative, one
concludes that $\det NS(S) = -48$.
\end{proof}

According to  the classification of Shioda-Inose \cite{shioda-inose},
there are four nonisomorphic singular $K3$ surfaces whose
N\'eron-Severi lattices have determinant $-48$.  The transcendental
lattice of such a surface is isomorphic to one of the lattices given
by the Gram matrices
\[
\begin{pmatrix}
2 & 0\\0 & 24
\end{pmatrix}, \qquad
\begin{pmatrix}
4 & 0\\0 & 12
\end{pmatrix}, \qquad
\begin{pmatrix}
6 & 0\\0 & 8
\end{pmatrix}, \qquad
\begin{pmatrix}
8 & 4\\4 & 8
\end{pmatrix}.
\]

\begin{prop}
The surface $S$ is isomorphic over $\Q(i,\w)$ to the quartic
surface in $\P^{3}$ defined by the equation $S':X'(X'{}^{3}+Y'{}^{3}) =
Z'(Z'{}^{3}+W'{}^{3})$.  As a
consequence, $S$ corresponds to the matrix $\begin{pmatrix}
8 & 4\\4 & 8 \end{pmatrix}$ in the Shioda-Inose classification.
\end{prop}

\begin{proof}
Consider the family of planes passing through the line $\ell_{5}$.
The intersection of $S_{0}$ and a plane in this family is $\ell_{5}$
and a plane cubic curve.  This yields another elliptic fibration
$\e_{2}:S\to\P^{1}$.  Setting $t=W/(X-iY)$ and using a standard
algorithm we can convert the fiber at the generic point to the
Weierstrass form
\[
y^2= x^{3} - 432t^2(t-1)^2(t+1)^2(t^2+1)^2,
\]
where the transformation is given by
\begin{alignat*}6
&X &&= &\ &72\sqrt{3} t(t^{2}+1)^2, &\qquad
&Y &&= &  &3y-36\sqrt{-3}t(t^4-1),\\
&Z &&= &  &-6\sqrt{-3}(t^{2}+1)x, &
&W &&= &  &-t \bigl(3 i y-36\sqrt{3} t (t^2+3) (t^2+1)\bigr).
\end{alignat*}

This fibration has six fibers of type IV at $t=0,\pm1,\pm i, \infty$.
The same fibration can be obtained by starting from
the surface $S'$ given by 
$$X'(X'{}^{3}+Y'{}^{3})= Z'(Z'{}^{3}+W'{}^{3}).$$ 
Namely, with $t=Z'/X'$, using the
transformation
\begin{alignat*}3
&X'&&= &\  & - 2 (t^4 + 1) y - t x^2 +12 t^3 x - 72 t (t^4-1)^2 , \\
&Y'&&= &  & - 2 (2 t^4 - 1) y + t x^2 -12 t^3 (t^4-1) x + 72 t (t^4-1)
,\\
&Z'&&= &  &\ t\bigl(- 2 (t^4 + 1) y - t x^2 + 12 t^3 x  - 72 t (t^4-1)^2
\bigr),\\
&W'&&= &  &\ 2 t (t^4 - 2) y + t^2 x^2  + 12 (t^4-1) x - 72 t^6 (t^4-1),
\end{alignat*} 
one finds exactly the same Weierstrass equation.  Since $S'$ is nonsingular
and minimal, we conclude that $S$ and $S'$ are isomorphic.

A result of Inose \cite{inose:quartic} shows that $S'$
corresponds to the matrix $\begin{pmatrix} 8 & 4\\4 & 8 \end{pmatrix}$
in the Shioda-Inose classification. This proves the assertion.
\end{proof}

\begin{remark}
The Shioda-Tate formula shows that the Mordell-Weil rank of the
fibration $\e_2$ is $20-2-6\times 2=6$.  However, since the
transformation from the $(X,Y,Z,W)$-coordinates to the affine
$(x,y,t)$-coordinates used here is not defined over $\Q$, one cannot
expect sections of $\e_2$ to correspond to solutions to our original
diophantine problem.
\end{remark}

\section{A polynomial solution and a base change}

{}From Proposition~\ref{MWL} one deduces that $E_t(\Q(t))$ is
generated by $\s_1$. A simple calculation shows that
\[
[2]\s_1 = \left(\frac{1}{4}t^2(t^2+8),
\frac{1}{8}t^2(t^4-20t^2-8)\right)
\]
is a polynomial solution to (\ref{weier}), but this does not
correspond to a polynomial solution to \eqref{eq_xyz}, as we
have
\[
(x,z) = \left(\frac{8(t^2-1)^2}{t^4-20t^2-8},
\frac{2t(t^2+8)(t^2-1)}{t^4-20t^2-8}\right).
\]
Since $[3]\s_1$, $[4]\s_1$, etc., do not yield a polynomial
solution to \eqref{eq_xyz}, we suspect that the elliptic fibration
$E_t$ gives no non-trivial parametric solutions to the
original question.  In fact, we can prove the following by an
elementary argument.
\begin{prop}\label{no-poly-solution}
The only polynomial solutions of the form $(f(t),t,g(t))$ to the
equation \eqref{eq_xyz} are $(\pm1,t,t)$.
\end{prop}

\begin{proof}
Suppose that a polynomial
$p(t)$ which is relatively prime to $t(t-1)(t+1)$ divides
$f(t)$.  Since $p(t)$ does not divide $t(f(t)^2+t^2-1)$,
its cube $p(t)^3$  divides $f(t)$.  Thus, we may write
$f(t)=t^\a(t-1)^\b(t+1)^\g f_0(t)^3$ ($0\le \a,\b,\g\le 2$).
Since $f_0(t)$ divides $g(t)$, we write $g(t)=f_0(t)g_0(t)$.
The polynomials $f_0(t)$ and $g_0(t)$ satisfy the equation
\[
t^{3\a+1}(t-1)^{3\b}(t+1)^{3\g}f_0(t)^6
+ t^{\a+1}(t-1)^{\b+1}(t+1)^{\g+1} = g_0(t)^3.
\]
The degree of the left hand side is
$\max\{3(\a+\b+\g)+6\deg f_0+1,\a+\b+\g+3\}$, but this
equals $3(\a+\b+\g)+6\deg f_0+1$ except in the case
$\a=\b=\g=\deg f_0=0$.  Since the degree of the right
hand side is a multiple of 3, the degree of the left hand
side cannot equal $3(\a+\b+\g)+6\deg f_0+1$.  Thus, we
have $\a=\b=\g=\deg f_0=0$.  In other words $f_0(t)$ must
be a constant and the only possibilities for this constant are
$\pm1$, which lead to the solutions $(\pm1,t,t)$.
\end{proof}

Observe that if we set $t=u^3$ in (\ref{weier}), then we find
a solution given as $(x_{1},y_{1})=(u^4(u^6-1),0)$.
This is a
2-torsion section of the elliptic surface corresponding to
\[
y_{1}^2=x_{1}^3-u^{12}(u^6-1)^3.
\]
Setting $y_{2}=y_{1}/u^6$ and $x_{2}=x_{1}/u^4$, a minimal
Weierstrass equation
\[
E'_u : y_{2}^2=x_{2}^3-(u^6-1)^3
\]
is obtained.
Let $\t$ be the above torsion section and let
$\s'_1$ be the section of $E'_u$ coming from $\s_1$.
We have that
\begin{align*}
\t &=  (u^6-1,0)\\
\s'_1 &= \bigl(u^2(u^6-1),(u^6-1)^2\bigr),
\end{align*}
and
\[
\pm\s'_1+\t = \bigl((u^2+2)(u^4+u^2+1),\mp3(u^4+u^2+1)^2\bigr).
\]
The latter sections correspond to
\[
\left\{\begin{aligned}
x&=\pm\frac{1}{3}(u^2-1)^2 \\[4pt]
y&=u^3 \\
z&=\pm\frac{1}{3}u(u^2-1)(u^2+2)
\end{aligned}\right.
\]
as solutions to the equation \eqref{eq_xyz}.
Taking the positive sign, and interchanging $x$ and $y$, 
one obtains the following solution to \eqref{eq_mkl}:
\[
\left\{\begin{aligned}
m&=\frac{1}{6}(u-1)(u^3-2u^2-4u-4)\\[4pt]
k&=u^{3} \\[4pt]
l&=\frac{1}{6}u(u^2-1)(u^2+2).
\end{aligned}\right.
\]
For integer values of $u$ not divisible by $3$, the 
$m$, $k$ and $l$
given above become integers. Thus, possibly after using the
symmetries mentioned in section~\ref{symm}, they give (infinitely
many nontrivial) solutions to the original diophantine problem.
For instance, $u=2$ leads to the solution
$3^3+4^3+5^3=6^3$. It should be remarked here that this
parametric solutions is in fact very old: according to
\cite[p.~582]{Di} it already appeared in a paper by
C. Pagliani which was published in 1829-1830.

Using a computer we found that \eqref{eq_xyz} has precisely
$32$ solutions in integers $0<y\leq x\leq 10^6$. Of these,
$15$ (almost half of them!) correspond to values of $u$ as
above.

\section{Further properties of $S$}

We will now study the geometry and arithmetic of the surface
$S$ a bit more closely. To this end, the elliptic fibration $E_t$
and the base change $E'_u$ will be used.
Denote by $E^-$ the elliptic curve over $\Q$ corresponding to the
Weierstrass equation $y^2=x^3-1$. Also, define a curve $C$ over $\Q$
by $s^2=u^6-1$. Note that $C$ is hyperelliptic and of genus $2$.
The hyperelliptic involution on $C$ will be written as $\io$; by
definition $\io(u,s)=(u,-s)$. Using the map $[-1]$ on $E^-$ as well,
one obtains an involution $\<\io\times[-1]\>$ on the product $C\times
E^-$.

\begin{prop}\label{quot1}
With notations as introduced above, the quotient
$(C\times E^-)/\<\io\times[-1]\>$ is birationally isomorphic over
$\Q$ to $E_u'$. 
\end{prop}

\begin{proof}
The function field of
$C\times E^-$ over $\Q$ is $\Q(u,x)[s,y]$ (with the relations
$y^2=x^3-1$ and $s^2=u^6-1$). The function field of the quotient
by $\io\times [-1]$ is the subfield of all functions invariant under
$(u,x,s,y)\mapsto(u,x,-s,-y)$. Over $\Q$, these invariants are
generated by $u, X'=s^2x$ and $Y'=s^3y$. The relation between
these is $Y'{}^2=X'{}^3-(u^6-1)^3$, so one finds precisely the
function field of $E_u'$ over $\Q$.
\end{proof}

Using a primitive cube root of unity $\w\in\overline{\Q}$
one defines $[\w]\in\mbox{Aut}(E^-)$ by 
$[\w](x,y)=(\w x,y)$. Similarly we take $\ph\in \mbox{Aut}(C)$
given by $\ph(u,s)=(\w^2u,s)$. Note that $[\w]$ and $\ph$ are
not defined over $\Q$, but the finite groups of automorphisms
they generate is. As a consequence, a quotient
such as $(C\times E^-)/\<\io\times[-1],\ph\times[\w]\>$ is
defined over $\Q$.

\begin{prop}
The surface $S$, or equivalently its elliptic
fibration $E_t$, is birationally isomorphic over
$\Q$ to $(C\times E^-)/\<\io\times[-1],\ph\times[\w]\>$.
\end{prop}

\begin{proof}
Consider
the surjection $E_u'\rightarrow E_t$ given by
\[
(u,X',Y')\mapsto (t,X,Y) = (u^3,u^4X',u^6Y').
\]
This is in fact the quotient map for the group of automorphisms
on $E_u'$ generated by $(u,X',Y')\mapsto(\w^2u,\w X',Y')$.
This generator can be lifted to the automorphism $\ph\times[\w]$ on 
$C\times E^-$. Using proposition~\ref{quot1} above now
finishes the proof.
\end{proof}

In the next section it will be explained how the above
description allows one to compute
the number of $\F_{p^n}$-rational points on the reduction
$S\bmod p$ in an easy way. As a preparation to that, we
now consider the second cohomology group $H^2(S)=H^2(S,\C)$.
The cycle class map allows one to view the N\'{e}ron-Severi group of
$S$ as a part of $H^2(S$; it generates a linear subspace
of dimension $20$ in our case. The orthogonal complement (with respect
to cup product) of this will be denoted $H^2_{\tr}(\tS)$. In terms of the
Hodge decomposition of $H^2(\tS)$, since $\tS$ is a singular $K3$ 
surface, one has
$H^2_{\tr}(\tS)=H^{2,0}(\tS)\oplus H^{0,2}(\tS)$ in the present
situation.
Using the relation between $\tS$ and $C\times E^-$ explained above,
one regards $H^2_{\tr}(\tS)$ as the subspace of 
$H^2_{\tr}(C\times E^-)$
on which $\io\times[-1]$ and $\ph\times[\w]$ act trivially.
Since the K\"{u}nneth components $H^2(C)\otimes H^0(E^-)$ and
$H^0(C)\otimes H^2(E^-)$ of $H^2(C\times E^-)$ are algebraic, this 
means that $H^2_{\tr}(\tS)$ can be seen as a subspace of $H^1(C)\otimes
H^1(E^-)$.

Note that the latter tensor product has dimension $8$. To understand it
even better we introduce the elliptic curve $E^+$, defined over $\Q$ by the
equation $\eta^2=\xi^3+1$. One has a morphism
$\pi_1\;:\;C\rightarrow E^+$, given by
$\pi_1(u,s)=(\xi=-u^{-2},\eta=su^{-3})$.
The global differential given by $\frac{\displaystyle
d\xi}{\displaystyle \eta}$
is pulled back to $2\frac{\displaystyle du}{\displaystyle s}$ under
$\pi_1$.
Hence it follows that under the pull back map $\pi_1^*$, the space
$H^1(E^+)$ is mapped to the $-1$-eigenspace in $H^1(C)$ for the
automorphism $\psi$ given by $\psi(u,s)=(-u,s)$.
Similarly, using $\pi_2\;:\;C\rightarrow E^-$ given by
$\pi_2(u,s)=(x=u^2,y=s)$ one computes
$\pi_2^*\frac{\displaystyle dx}{\displaystyle y}=2u\frac{\displaystyle
du}
{\displaystyle s}$ and concludes that $\pi_2^*$ maps $H^1(E^-)$ to the
$+1$-eigenspace in $H^1(C)$ for $\psi$.
Hence using $(\pi_1\times\pi_2)^*$ one identifies
$$H^1(C)\otimes H^1(E^-)=\left(H^1(E^+)\otimes H^1(E^-)\right)\oplus
\left(H^1(E^-)\otimes H^1(E^-)\right).$$

We claim that under this identification, the subspace $H^2_{\tr}(\tS)$
coincides with the transcendental part $H^2_{\tr}(E^+\times E^-)$
inside $H^1(E^+)\otimes H^1(E^-)$. To see this, note that the latter
transcendental part is generated by the holomorphic form
$\frac{\displaystyle d\xi}{\displaystyle \eta}\otimes
\frac{\displaystyle dx}{\displaystyle y}$ and its anti-holomorphic
complex conjugate. Using $\pi_1^*$, this holomorphic form corresponds to
$2\frac{\displaystyle du}{\displaystyle s}\otimes
\frac{\displaystyle dx}{\displaystyle y}$, which is invariant under
both $\io\times[-1]$ and $\ph\times[\w]$, hence it (and its complex
conjugate)
is in $H^2_{\tr}(\tS)$. Since this space is $2$-dimensional, the claim
is proven.

\vspace{\baselineskip}\noindent
{\sl Remark.} The proofs provided in this section in fact describe
an inclusion of function fields $\Q(C\times E^-)\supset\Q(E'_u)\supset
\Q(E_t)$. It is easily verified that the composition of these
inclusions corresponds to the finite (in fact, cyclic of degree $6$)
map $\psi:C\times E^- \rightarrow E_t=S$ which assigns to a point
$(u,s),\,(x_0,y_0)\in C\times E^-$ the point with coordinates 
$x=s/y_0, y=u^3, z=usx_0/y_0$ on $S$. The cyclic group of order $6$
for which we take the quotient here, turns out to have precisely $16$ 
orbits in $C\times E^-$ consisting of less than $6$ points. Blowing
up the quotient singularities obtained in this way, yields $16$
(independent) elements in the N\'eron-Severi group of $S$. In fact
this provides an alternative way to see that $S$ has N\'eron-Severi rank 
$20$: One can lift the $4$ generators of the N\'eron-Severi group of
$E^-\times E^-$ to cycles in $C\times E^-$, and push them forward
to $S$. In this way $20$ independent elements are obtained, as
follows immediately by computing their intersection numbers.

As a final remark, consider the morphism $C\rightarrow E^-$ given
by $P\mapsto \pi_2(P)+(0,1)$. The graph of this morphism defines a
curve in $C\times E^-$, and the image in $E_t$ of this curve under $\psi$
is precisely the `polynomial section' described by Pagliani in 1829/30
and by us in the previous section.

\section{Counting points on $S$ over finite fields}
Since $S$ is a singular $K3$ surface, one knows from
\cite{shioda-inose} that for sufficiently large extensions $K$
of the finite prime field $\F_p$ one can describe the number
of $K$-rational points on the reduction of $S$ in terms of
a Hecke character. However, it is not easy to describe,
for which field extensions $K$ their result gives the number
of points. For the present example, we present a Hecke character
which gives the number of rational points over {\em all} finite
fields. Although this does not seem to have an immediate relation
to the problem of sums of consecutive cubes, it is included here
because it is an interesting and nontrivial problem for singular
$K3$'s over number fields in general, and not many interesting
examples where this has been done have been published.
One may note that the present case turns out to be
much simpler than, e.g., the example described in \cite{P-T-vdV}.

To compute the number of points on $S$ over finite fields, the
$\ell$-adic cohomology groups $H^i=H^i(S_{\overline{\Q}},\Q_{\ell})$
are used. These spaces are
$G_{\Q}=\mbox{Gal}(\overline{\Q}/\Q)$-modules.
One knows $H^0\cong\Q_{\ell}$, with the trivial Galois action and
$H^4\cong\Q_{\ell}(-2)$, which means that $H^4$ is a $\Q_{\ell}$-vector
space
of dimension $1$ and a Frobenius element $\s_p\in G_{\Q}$ at a 
prime $p$ acts on it by multiplication by $p^2$. Furthermore, since $S$ is an
elliptic surface with base $\P^1$ it follows that $H^1=H^3=(0)$.
Given a Frobenius element $\s_p\in G_{\Q}$ at a `good' prime $p\neq\ell$
(in our situation `good' means that one of the equations defining $S$
gives a $K3$ surface over $\F_p$ as well; this happens when $p\neq
2,3$),
the number of points over $\F_{p^n}$ on the reduction mod $p$ of $S$
is given by the Lefschetz trace formula
$$\sum_{i=0}^4
(-1)^i\mbox{trace}(\s_p^n|H^i(S_{\overline{\Q}},\Q_{\ell}))
=1+p^{2n}+\mbox{trace}(\s_p^n|H^2(S_{\overline{\Q}},\Q_{\ell})).$$
We will make this explicit by studying
$H^2(S_{\overline{\Q}},\Q_{\ell})$.

There exists a $G_{\Q}$-equivariant cycle class map which injects the
N\'{e}ron-Severi group of $S$ into
$H^2(S_{\overline{\Q}},\Q_{\ell}(1))$, which is the same space 
$H^2$ but with a
different Galois action: $\s_p$ acts on
$H^2(S_{\overline{\Q}},\Q_{\ell})$
exactly as $p\s_p$ does on $H^2(S_{\overline{\Q}},\Q_{\ell}(1))$.
Write $H^2_{\mbox{\scriptsize alg}}$ for the subspace of
$H^2(S_{\overline{\Q}},\Q_{\ell})$ generated by the
N\'{e}ron-Severi group, and $H^2_{\mbox{\scriptsize tr}}$ for the
orthogonal complement. Both spaces are $G_{\Q}$-invariant, and of course
$$\mbox{trace}(\s_p^n|H^2(S_{\overline{\Q}},\Q_{\ell}))=
\mbox{trace}(\s_p^n|H^2_{\mbox{\scriptsize alg}})+
\mbox{trace}(\s_p^n|H^2_{\mbox{\scriptsize tr}}).$$

\begin{prop}
For $p\neq 2,3,\ell$ prime and $\s_p\in G_{\Q}$ a Frobenius element at
$p$,
one has
\[
\mbox{trace}(\sigma_p^n|H^2_{\mbox{\scriptsize
alg}})=16p^n+3{\leg{-3}{p}}^{n}{p^n}
+{\leg{-4}{p}}^n{p^n}.
\]
\end{prop}

\begin{proof} By what is said above, it suffices  to compute the action
of a Frobenius element $\s_p$ on a set of generators of the
N\'{e}ron-Severi
group. Using the elliptic fibration $E_t$, we study such generators.

Firstly, there are the zero section and a fiber. On each of these 
$\s_p$ acts trivially, so they contribute $p^n+p^n=2p^n$ to our trace.

Next we consider the sections $\s_1$ and $[\w]\s_1$.
Note that since
$\s_1+[\w]\s_1+[\overline{\w}]\s_1=0$ in the group law, it follows that
$\s_p([\w]\s_1)=-\s_1-[\w]\s_1$ in case ${\leg{-3}{p}}=-1$. On the other
hand,
$\s_1$ is fixed under every element of $G_{\Q}$. So it follows that
these two sections contribute $\left(1+{\leg{-3}{p}}^n\right)p^n$ to 
the trace.

It remains to compute the contributions from the irreducible
components
of the singular fibers, which do not meet the zero section.
These fibers are over $t=0$, $t=\pm 1$ and $t=\infty$.
Over $t=0$, the fibre is of type $IV^*$ and all its
components turn out to be rational.
Hence they contribute $6p^n$ to our trace.

Over $t=\pm 1$ one finds fibers of type $I^*_0$. Of the components
not meeting the zero section, one is always rational. With $T$ a
coordinate on this one, the other three meet it in points satisfying 
$T^3-(\pm 2)^3=0$.
Hence two of them are interchanged by $\mbox{Gal}(\Q(\w)/\Q)$ and the
third one is rational. This means that one obtains
\[
2\left(p^n+p^n+\left(1+{\leg{-3}{p}}^n\right)p^n\right)
=6p^n+2{\leg{-3}{p}}^{n}{p^n}
\]
from these fibers to the trace.

To study the fiber over $t=\infty$, one changes coordinates using
${\tilde y}=t^{-6}y, {\tilde x}=t^{-4}x$ and $s=t^{-1}$.
The elliptic fibration $E_t$ is then given by ${\tilde y}^2={\tilde
x}^3-s^2(1-s^2)^3$.
Over $t=\infty$, which corresponds to $s=0$, one finds a fiber of type
IV.
The two components not meeting the zero section are interchanged by
$\mbox{Gal}(\Q(i)/\Q)$. Hence they add $(1+{\leg{-4}{p}}^n)p^n$ to the
trace.

Summing all contributions now proves the proposition.
\end{proof}

\vspace{\baselineskip}
The $G_{\Q}$-space $H^2_{\tr}$ will be our next object to study. Recall
that we
already showed it is related to $H^2_{\tr}(A_{\overline{\Q}},\Q_{\ell})$
in which $A=E^+\times E^-$. In fact,
this was only done for complex cohomology, but using a comparison
theorem and noting that the morphisms $C\times E^-\rightarrow S$ and
$C\rightarrow E^+\times E^-$ we used are defined over $\Q$, the same
holds for $\ell$-adic cohomology. So one concludes
$$\mbox{trace}(\s_p^n|H^2_{\tr})=
\mbox{trace}(\s_p^n|H^2_{\tr}(A_{\overline{\Q}},\Q_{\ell})).$$
It is relatively standard how this latter trace is computed, as will
be explained now.

It has been shown that $H^2_{\tr}\subset H^1(E^+)\otimes H^1(E^-)$.
This K\"{u}nneth component in $\ell$-adic cohomology is
$G_{\Q}$-invariant,
for instance because its orthogonal complement, which is generated by
the cycle classes of $\{0\}\times E^-$ and $E^+\times\{0\}$, obviously
is.
Note that $\Z[\w]$ is the endomorphism ring of both $E^+,E^-$, and all
these
endomorphisms are defined over $\Q(\w)$. Hence the $2$-dimensional
$\Q_{\ell}$-vector
spaces $H^1(E^{\pm})$ are in fact free rank $1$ modules over
$\Q_{\ell}\otimes\Z[\w]$. In particular, this means that when we
restrict
the Galois action to $G_{\Q(\w)}=\mbox{Gal}(\overline{\Q}/\Q(\w))$, then
it
becomes abelian, i.e., it factors over the Galois group of the maximal
abelian extension of $\Q(\w)$. By class field theory this implies that
our
actions on $H^1(E^{\pm})$ are given by Hecke characters, which will be
denoted
$\chi^+,\chi^-$, respectively.
These characters will now be described.

Write $K=\Q(\w)$. The Hecke characters $\chi^{\pm}$ can be thought of as
products $\prod_v \chi_v^{\pm}$ in which the product is taken over all
places of $K$ (including the infinite one). Here $\chi_v$ is a
multiplicative character; for $v=\infty$ it is given by
$$\chi_{\infty}^{\pm}\;:\;\C^*\longrightarrow
\C^*\;\;:\;\;z\mapsto\frac1z.$$
For finite places one has $\chi_v^{\pm}\;:\;K_v^*\rightarrow\Q(\w)^*$.
These characters have the property that for a place $v$ not dividing
$(2)$ or $(3)$, and $v$ corresponding to a (principal) prime ideal $(\pi)$ of
norm $q$ in $\Z[\w]$,
one has that $\chi_v^{\pm}$ sends any uniformizing element of $K_v$ to
the generator of $(\pi)$ which as an element of the endomorphism ring
of $E^{\pm}$ gives after reduction modulo $(\pi)$ the $q$th power
endomorphism. In particular, for such $v$ one has that $\chi_v^{\pm}$
is $1$ on all elements of $K_v$ which have valuation $0$.
If $\s_p\in G_{\Q}$ is a Frobenius element at a prime $p\neq 2,3$,
then (for $n$ non-zero)
$$\mbox{trace}(\s_p^{f_pn}|H^1(E^{\pm}))=\sum_{v|p}\chi_v^{\pm}(\pi_v)^n
,$$
in which $\pi_v\in K_v$ is any uniformizing element and where $f_p$ is the
degree of the extension $K_v/\Q_p$. Moreover,
$\mbox{trace}(\s_p^m||H^1(E^{\pm})) =0$ when $f_p$ does not divide $m$.

An explicit description of the $\chi_v^{\pm}$ now boils down to finding
the $q$th power endomorphism on $E^{\pm}$ over $\F_q$. In case $p\neq 2$
is a prime number $\equiv 2\bmod 3$, the square of the $p$th power map has to
be considered, and this equals $[-p]$ as endomorphism (the reduction of
both $E^+,E^-$ is supersingular in this case).
When $p\equiv 1\bmod 3$, both primes above $p$ in $\Z[\w]$ are primes of
ordinary reduction for $E^{\pm}$. So here the full endomorphism ring in
characteristic $p$ is $\Z[\w]$, and we have to find out which
element the $p$th power map corresponds to. Firstly, this map has degree
$p$, so we look for an element of norm $p$. Next, the map is inseparable, so
we want an element whose reduction modulo the prime of $\Z[\w]$ under
consideration is $0$. This implies we want a generator $\pi$ of the
prime ideal under consideration. To find out which generator, consider some 
torsion points on $E^{\pm}$. Note that since $p\equiv 1\bmod 3$, all the
$2$-torsion points on both $E^{\pm}$ are rational over $\F_p$. 
Hence the $p$th power map $\pi$ fixes the $2$-torsion, which implies that 
$2|(\pi-1)$. Moreover, both $[\w]$ and $[\overline{\w}]$ act the same way 
(in fact, act trivially) on the points with $x-$ or $\xi$-coordinate $0$. 
Hence these points are in the kernel of $[\w-\overline{\w}]=[\sqrt{-3}]$. 
This map has degree $3$, so we found all the $[\sqrt{-3}]$-torsion. 
On $E^+$, these points are $\F_p$-rational, so it follows 
$\pi\equiv 1\bmod 2\sqrt{-3}$ in this case. This determines $\pi$.
For $E^-$, the $p$th power map acts as $[+1]$ on the
$[\sqrt{-3}]$-torsion when $p\equiv 1\bmod 12$, so in this case again 
$\pi$ is determined by $\pi\equiv 1\bmod 2\sqrt{-3}$. 
However, when $p\equiv 7\bmod 12$ then the $p$th power map acts as $[-1]$ 
on the $[\sqrt{-3}]$-torsion, hence $\pi$ is determined by 
$\pi\equiv -1\bmod 2\sqrt{-3}$.

For completeness, and to allow us to relate $H^2_{\tr}$ to a modular
form later on, we also describe the $\chi^{\pm}_v$ for $v|(2)(3)$.
These can be found using that $\prod_v \chi_v(x)=1$, where $x\in \Q(\w)$
and where the product is over all places $v$ of $\Q(\w)$.
Furthermore, the local units $\Z_2[\w]^*$ and $\Z_3[\w]^*$ are known to
have finite image, hence we know that this image has to be inside the
$6$th roots of unity. So in particular all sufficiently small subgroups
$1+2^n\Z_2[\w]$ resp. $1+(\sqrt{-3})^m\Z_3[\w]$, which consist of only
$6$th powers, are mapped to $1$. Since there is only one prime
above each of $(2)$, $(3)$, we denote the associated local characters
by $\chi_2^{\pm}$ and $\chi_3^{\pm}$ respectively. One computes:
\begin{enumerate}
\item For $p=3$ one identifies
\[
\Z_3[\w]^*/1+(\sqrt{-3})\cong\left(\Z_3[\w]/\sqrt{-3}\right)^*
\cong(\Z/3\Z)^*\cong\{\pm 1\}.
\]
Then both $\chi_3^+$ and $\chi_3^-$ are trivial on $1+(\sqrt{-3})$, and
on $\Z_3[\w]^*$ they are given using the identification above.
Next, $\chi_3^{\pm}(\sqrt{-3})=\sqrt{-3}$.
\item For $p=2$, identify
$$\Z_2[\w]^*/1+(2)\cong
\left(\Z_2[\w]/2\right)^*\cong\{1,\w,\overline{\w}\}.$$
Then $\chi_2^+$ is trivial on $1+2\Z_2[\w]$, and it is given by this
identification on $\Z_2[\w]^*$. Moreover, $\chi_2^+(2)=-2$.

The character $\chi_2^-$ turns out to be trivial on $1+4\Z_2[\w]$, and on
$\Z_2[\w]^*$ it is given by
$$\chi_2^-(u)=(-1)^{(u\overline{u}-1)/2}\chi_2^+(u).$$
Finally, $\chi_2^-(2)=-2$.
\end{enumerate}

\vspace{\baselineskip}
{}From the discussion above it is clear that the restriction to
$G_{\Q(\w)}$ of the $G_{\Q}$-representation $H^2_{\tr}(A)$ is given by
the product $\chi=\chi^+\chi^-$. Namely, over $\Q(\w,i)$ one finds
two copies of  $G_{\Q(\w,i)}$-representation in
$H^1(E^+)\otimes H^1(E^-)$, corresponding to the graph of
an isomorphism between these elliptic curves, and the composition of that
graph with $[\w]$. These representations correspond to the Hecke
character (in fact, Dirichlet character) $\chi^+\overline{\chi^-}$.
Hence the remaining transcendental part corresponds to $\chi$.

The explicit description of $\chi^{\pm}$ shows that $\chi_3$ is in fact
trivial on the units $\Z_3[\w]^*$. In fact, one computes that $\chi$
is a Hecke character of conductor $(2)^4$. Since such Hecke characters
are well known to be intimately related to modular forms,
one concludes using e.g. \cite[Thm. 2.4.2]{top-thesis} That the
$L$-series attached to $H^2_{\tr}$, which equals
\[
L(s,\chi)=\prod_{\begin{subarray}{c} v\; \text{finite}\\  v\neq
2\end{subarray}}
(1-\chi_v(\pi_v)N\pi_v^{-s})^{-1},
\]
is in fact the $L$-series of a cusp form of weight $3$ and character
${\leg{-3}{\cdot}}$ for $\Gamma_0(48)$.

If one works out the correspondence between this cusp form $f$ and the
$L$-series in more detail, one finds that $f$ is given by the
$q$-expansion
\[
\frac{1}{6}\sum_{
\begin{subarray}{c}
m,n\in {\bf Z}\\ (m,n)\not\equiv(0,0)\bmod 2
\end{subarray}}
(m+n\w)^2\chi_2(n+m\w)^{-1}q^{m^2-mn+n^2}.
\]
Hence one finds
\begin{multline*}
f=q+3q^3-2q^7+9q^9-22q^{13}-26q^{19}-6q^{21}+25q^{25}\\
+27q^{27}+46q^{31}+26q^{37}-66q^{39}+22q^{43}-45q^{49}\ldots
\end{multline*}

It was pointed out to us by Ken Ono that this cusp form
has the following description.
Denote by $\eta(z)$ the Dedekind eta-function, i.e., the function
given by $\eta(z)=e^{2\pi i z/24}\prod_{n\geq 1}(1-e^{2\pi i n z})$.
Using the famous transformation rules for $\eta$ it is not hard
to verify that 
$$g(z):= \frac{\displaystyle \eta(12 z)^9\eta(4z)^9}
{\displaystyle \eta(2z)^3\eta(6z)^3\eta(8z)^3\eta(24z)^3}
$$
is a modular form of weight $3$ and character
${\leg{-3}{\cdot}}$ for $\Gamma_0(48)$ as well. Comparing
coefficients in fact shows that $f=g$, hence our $L$-series is the
one associated with this product of eta-functions.

\vspace{\baselineskip}
The above discussion allows one to compute the number of points on
the somewhat abstractly defined variety $S$ over finite fields.
However, since it is a purely combinatorial matter
(compare\cite{vangeeman-top}) to describe the
difference between $S(\F_{p^n})$ and the set
$\{(t,x,y)\in\F_{p^n}^{(3)}\;|\; y^2=x^3-t^4(t^2-1)^3\}$ we will for
explicitness state the final result of the above discussion in terms of
this explicit equation.
\begin{theorem}
Suppose $p\geq 5$ is a prime number and $n>0$ an integer.
The number $N(p,n)$ of solutions $(t,x,y)\in\F_{p^n}^{(3)}$ to the
equation $y^2=x^3-t^4(t^2-1)^3\}$ is given by
$$N(p,n)=p^{2n}+p^n+{\leg{-3}{p}}^np^n+a_{p^n},$$
in which $a_{p^n}$ is given by any of the two following descriptions.
\begin{enumerate}
\item $a_{p^n}$ is the coefficient of $q^{p^n}$ in the cusp form $f$
of weight $3$ and character ${\leg{-3}{\cdot}}$ for $\Gamma_0(48)$
given by
$$q \prod_{n\geq 1} \frac{\displaystyle (1-q^{12n})^9(1-q^{4n})^9}
{\displaystyle (1-q^{2n})^{3}
(1-q^{6n})^{3}(1-q^{8n})^{3}(1-q^{24n})^{3}};$$
\item For $p\equiv 2\bmod 3$ one has $a_{p^n}=0$ whenever $n$ is odd,
and $a_{p^n}=p^{n}$ whenever $n$ is even.

For $p\equiv 1\bmod 3$, write $p=m^2-mn+n^2$. Then
$a_{p^n}=\alpha^n+\overline{\alpha}^n$, in which
$\alpha={\leg{-4}{p}}\w^a(m+n\w)^2$, and where the exponent $a$ is chosen
such that $m+n\w-w^a\in 2\Z[\w]$.
\end{enumerate}

\end{theorem}

\providecommand{\bysame}{\leavevmode\hbox to3em{\hrulefill}\thinspace}


\begin{thebibliography}{P-T-vdV}

\bibitem[Di]{Di}
L.E. Dickson, History of the theory of numbers, Volume II.
Washington: Carnegie Institution, 1920.

\bibitem[Ino]{inose:quartic}
H.~Inose, \emph{On certain {Kummer} surfaces which can be realized as
  non-singular quartic surfaces in {$\P^3$}}, J. Fac. Sci. Univ. Tokyo
    \textbf{23} (1976), 545--560.
    
\bibitem[Kuw]{kuwata:can-height}
M.~Kuwata, \emph{Canonical height and elliptic {$K$3} surfaces}, J.
Number Theory \textbf{36} (1990), 399--406.

\bibitem[K-T]{K-T}
M. Kuwata and J. Top, \emph{An elliptic surface related to sums of
consecutive squares},
Expos. Math. \textbf{12} (1994), 181--192.

\bibitem[M-P]{miranda-persson}
R.~Miranda and U.~Persson, \emph{Torsion subgroup of elliptic surfaces},
Compositio Math. \textbf{72} (1989), 249--267.

\bibitem[P-T-vdV]{P-T-vdV}
C.~Peters, J.~Top, and M.~van~der Vlugt, \emph{The Hasse zeta function
of a  {$K3$} surface related to the number of words of weight $5$ in the
{Melas} codes}, J. Reine Angew. Math. \textbf{432} (1992), 151--176.

\bibitem[P-R]{P-R}
S.~Platiel and J.~Rung, \emph{Nat\"{u}rliche Zahlen als Summen
aufeinander folgender Quadratzahlen}, Expos. Math. \textbf{12}
(1994), 353--361.

\bibitem[Shi]{shioda:MWL}
T.~Shioda, \emph{On the {Mordell-Weil} lattices}, Comment. Math. Univ.
Sancti Pauli \textbf{39} (1990), 211--240.

\bibitem[S-I]{shioda-inose}
T.~Shioda and H.~Inose, \emph{On singular {$K$3} surfaces}, Complex Analysis
  and Algebraic Geometry (1977), 119--136.
  
\bibitem[Str]{Str}
R.J.~Stroeker, \emph{On the sum of consecutive cubes being a
perfect square}, Compositio Math. \textbf{97} (1995), 295--307.

\bibitem[Top]{top-thesis}
J.~Top, \emph{Hecke {$L$-series} related with algebraic cycles or with
{Siegel} modular forms}, Ph.D. thesis, University of Utrecht, 1989.

\bibitem[vG-T]{vangeeman-top}
B.~van Geemen and J.~Top, \emph{Selfdual and non-selfdual
3-dimensional Galois representations},
Compositio Math. \textbf{97} (1995), 51--70.
\end{thebibliography}
\end{document}